\documentclass{article}
\usepackage{amsmath,amsfonts,amssymb,bm,hyperref}
\numberwithin{equation}{section}
\newtheorem{theorem}{Theorem}[section]
\newtheorem{lemma}{Lemma}[section]
\newtheorem{corollary}{Corollary}[section]

\newtheorem{remark}{Remark}[section]

\usepackage[T1]{fontenc}
\usepackage[utf8]{inputenc}
\usepackage{authblk}

\author{Zhirayr Avetisyan\thanks{Z.Avetisyan@ucl.ac.uk;
Zhirayr Avetisyan was supported by EPSRC grant EP/M000079/1}}
\author{Yan-Long Fang\thanks{Yan.Fang.12@ucl.ac.uk}}
\author{Dmitri Vassiliev\thanks{D.Vassiliev@ucl.ac.uk,
\url{http://www.homepages.ucl.ac.uk/\~ucahdva/};
Dmitri Vassiliev was supported by EPSRC grant EP/M000079/1}}
\affil{Department of Mathematics,
University College London,\\ Gower Street, London WC1E 6BT, UK}

\begin{document}

\title{Spectral asymptotics for first order systems}

\maketitle
\begin{abstract}
This is a review paper outlining recent progress in the spectral analysis of first order systems.
We work on a closed manifold and study an elliptic self-adjoint first order
system of linear partial differential equations.
The aim is to examine the spectrum and derive asymptotic formulae
for the two counting functions.
Here the two counting functions are those for the positive and the negative eigenvalues.
One has to deal with positive and negative eigenvalues separately
because the spectrum is, generically, asymmetric.
\end{abstract}

\section{Introduction}
\label{Introduction}

This paper was inspired by Yuri Safarov's treatment of first order systems in \cite{SafarovDSc}.
Safarov was one of the first researchers to attempt a detailed spectral analysis of
first order systems on closed manifolds. A brief historical review can be found
in Section~11 of \cite{jst_part_a}.

Let $L$ be a first order linear differential operator acting on $m$-columns
of complex-valued half-densities over a connected
closed (i.e.~compact and without boundary) $n$-dimensional  manifold $M$.
Throughout this paper we assume that $m,n\ge2$.

Let $L_\mathrm{prin}(x,p)$ and $L_\mathrm{sub}(x)$ be the principal
and subprincipal symbols of $L$.
Here $x=(x^1,\ldots,x^n)$ denotes local coordinates
and $p=(p_1,\ldots,p_n)$  denotes the dual variable (momentum).
The principal and subprincipal symbols
are defined in the same way as for scalar operators,
see subsection 2.1.3 in \cite{mybook},
only now they are $m\times m$ matrix-functions
on $T^*M$ and $M$ respectively.
As our operator $L$ is first order and differential (as opposed to pseudodifferential),
the principal and subprincipal symbols uniquely determine the operator.
In other words, the principal and subprincipal symbols
provide an invariant analytic description of our first
order differential operator $L$.

We assume our operator $L$ to be formally self-adjoint (symmetric)
with respect to the standard inner product on $m$-columns
of complex-valued half-densities,
which is equivalent to the
assumption that the principal and subprincipal symbols are Hermitian.
We also assume that our operator $L$ is elliptic:
\begin{equation}
\label{definition of ellipticity}
\det L_\mathrm{prin}(x,p)\ne0,\qquad\forall(x,p)\in T^*M\setminus\{0\}.
\end{equation}
Condition \eqref{definition of ellipticity}
and the fact that dimension $n$ is greater than or equal
to two imply that $m$, the number of equations, is even.
See Remark 2.1 in \cite{israel} for details.

Under the above assumptions $L$ is a self-adjoint operator in
$L^2(M;\mathbb{C}^m)$ with domain $H^1(M;\mathbb{C}^m)$ and
the spectrum of $L$ is discrete, accumulating to
$+\infty$ and to $-\infty$.
Let $\lambda_k$ and $v_k(x)$ be
the eigenvalues and the orthonormal eigenfunctions of the operator~$L$; the
particular enumeration of these eigenvalues
(accounting for multiplicities)
is irrelevant for our purposes.
Each $v_k(x)$ is, of course,
an $m$-column of half-densities.

The main objective of this paper is to derive asymptotic formulae for the distribution
of large, in terms of modulus, eigenvalues of $L$. We will deal with positive and
negative eigenvalues separately and we will see that the asymptotic distribution of positive
and negative eigenvalues is somewhat different.
Note that asymmetry of the spectrum (with respect to zero) is a major subject
in geometry, see, for example,
\cite{atiyah_short_paper,atiyah_part_1,atiyah_part_2,atiyah_part_3}.

\section{The propagator}
\label{The propagator}

In this paper, as in \cite{mybook}, we use the wave equation method
(Levitan's method) to derive spectral asymptotics.
In our case, when we are dealing with an operator which is not semi-bounded,
the wave equation method is especially natural.
Furthermore, the wave equation method is the only physically meaningful way
of introducing a time coordinate when dealing with a first order elliptic system.

Let $h^{(j)}(x,p)$ be the eigenvalues
of the matrix-function
$L_\mathrm{prin}(x,p)$.
Throughout this paper we assume that these are
simple for all $(x,p)\in T^*M\setminus\{0\}$.
We enumerate the $h^{(j)}(x,p)$ in increasing order, using
a negative index $j=-m/2,\ldots,-1$ for negative $h^{(j)}(x,p)$
and
a positive index $j=1,\ldots,m/2$ for positive $h^{(j)}(x,p)$.
By $v^{(j)}(x,p)$ we denote the corresponding normalised eigenvectors.
Note that as our operator $L$ is first order and differential (as opposed to pseudodifferential) we have the following symmetry:
\begin{equation}
\label{symmetry property 1}
h^{(-j)}(x,p)=-h^{(j)}(x,-p),
\quad
v^{(-j)}(x,p)=v^{(j)}(x,-p),
\quad
j=1,\ldots,m/2.
\end{equation}

Now, let $x^{n+1}\in\mathbb{R}$ be the additional `time' coordinate.
Consider the Cauchy problem
\begin{equation}
\label{initial condition most basic}
\left.w\right|_{x^{n+1}=0}=v
\end{equation}
for the hyperbolic system
\begin{equation}
\label{dynamic equation most basic}
(-i\partial/\partial x^{n+1}+L)w=0
\end{equation}
on $M\times\mathbb{R}$.
The $m$-column of half-densities $v=v(x^1,\ldots,x^n)$ is given
and the $m$-column of half-densities $w=w(x^1,\ldots,x^n,x^{n+1})$ is
to be found.
The solution of the Cauchy problem
\eqref{initial condition most basic},
\eqref{dynamic equation most basic}
can be written as $\,w=U(x^{n+1})\,v$, where
\begin{multline}
\label{definition of wave group}
U(x^{n+1}):=e^{-ix^{n+1}L}
\\
=\sum_k e^{-ix^{n+1}\lambda_k}\,v_k(x^1,\ldots,x^n)\int_M[v_k(y^1,\ldots,y^n)]^*(\,\cdot\,)\,dy^1\ldots dy^n
\end{multline}
is the \emph{propagator}; here and further on the star stands for Hermitian conjugation.
The propagator $U(x^{n+1})$ is a one-parameter
family of unitary operators.

\begin{remark}
\label{remark on x n plus 1}
We chose to denote the `time' coordinate by $x^{n+1}$ rather than
by $t$ because some constructions presented in the current paper
work in the relativistic setting, i.e.~when there is no distinguished time
direction and the coordinates $x^1,\ldots,x^n$ and $x^{n+1}$ are
`mixed up'. Such an approach was pursued in \cite{nongeometric}
and, to a certain extent, in \cite{israel}.
\end{remark}

It was shown by Safarov \cite{SafarovDSc}
that the propagator can be constructed explicitly,
modulo $C^\infty$,
as a sum of $m$ oscillatory integrals
(Fourier integral operators)
\begin{equation}
\label{wave group as a sum of oscillatory integrals}
U(x^{n+1})\overset{\operatorname{mod}C^\infty}=
\sum_j
U^{(j)}(x^{n+1})\,,
\end{equation}
where the phase function of each oscillatory integral
$U^{(j)}(x^{n+1})$ is associated with the corresponding
Hamiltonian $h^{(j)}(x^1,\ldots,x^n,p_1,\ldots,p_n)$
and summation is carried out over all
nonzero integers $j$ from $-m/2$ to $+m/2$.
The notion of a phase function associated with a Hamiltonian is
defined in Section 2 of \cite{jst_part_a}
and Section~2.4 of \cite{mybook}.
Safarov's initial exposition \cite{SafarovDSc} of the construction leading up to
\eqref{wave group as a sum of oscillatory integrals} was quite concise.
A more detailed exposition was later
given in \cite{NicollPhD} and \cite{jst_part_a}.

We will now state the two main results regarding the properties
of the oscillatory integrals $U^{(j)}(x^{n+1})$ appearing in the RHS of formula
\eqref{wave group as a sum of oscillatory integrals}.
From this point till the end of the section we assume the index $j$ to be fixed.

The first result concerns the principal symbol of
the oscillatory integral $U^{(j)}(x^{n+1})$.
The notion of the principal symbol of an oscillatory
integral is defined in accordance with Definition 2.7.12 from
\cite{mybook}.
The principal symbol of the oscillatory
integral $U^{(j)}(x^{n+1})$ is a complex-valued
$m\times m$ matrix-function on
$\,M\times\mathbb{R}\times(T^*M\setminus\{0\})\,$.
We denote the arguments of this principal symbol by
$x^1,\ldots,x^n$ (local coordinates on~$M$),
$x^{n+1}$ (`time' coordinate on $\mathbb{R}$),
$y^1,\ldots,y^n$ (local coordinates on $M$)
and
$q_1,\ldots,q_n$ (variable dual to $y^1,\ldots,y^n$).

Further on in this section and the next section we use $x$, $y$, $p$ and $q$ as shorthand for
$\,x^1,\ldots,x^n$,
$\,y^1,\ldots,y^n$,
$\,p_1,\ldots,p_n$
and
$\,q_1,\ldots,q_n$
respectively.
The additional `time' coordinate $x^{n+1}$ will always be written separately.

In order to write down the principal symbol of the oscillatory
integral
\linebreak
$U^{(j)}(x^{n+1})$ we need to introduce some notation.

Curly brackets will denote the Poisson bracket on matrix-functions
$\{P,R\}:=P_{x^\alpha}R_{p_\alpha}-P_{p_\alpha}R_{x^\alpha}$
and its further generalisation
\begin{equation}
\label{definition of generalised Poisson bracket}
\{F,G,H\}:=F_{x^\alpha}GH_{p_\alpha}-F_{p_\alpha}GH_{x^\alpha},
\end{equation}
where the subscripts $x^\alpha$ and $p_\alpha$
indicate partial derivatives and
the repeated index $\alpha$ indicates summation over $\alpha=1,\ldots,n$.

We define the scalar function
$f^{(j)}:T^*M\setminus\{0\}\to\mathbb{R}$ in accordance with the formula
\begin{equation*}
\label{phase appearing in principal symbol}
f^{(j)}:=[v^{(j)}]^*L_\mathrm{sub}v^{(j)}
-\frac i2
\{
[v^{(j)}]^*,L_\mathrm{prin}-h^{(j)},v^{(j)}
\}
-i[v^{(j)}]^*\{v^{(j)},h^{(j)}\}.
\end{equation*}

By
$(x^{(j)}(x^{n+1};y,q),p^{(j)}(x^{n+1};y,q))$ we denote the Hamiltonian trajectory
originating from the point $(y,q)$, i.e.~solution of the system of
ordinary differential equations (the dot denotes differentiation in $x^{n+1}$)
\begin{equation*}
\label{Hamiltonian system of equations}
\dot x^{(j)}=h^{(j)}_p(x^{(j)},p^{(j)}),
\qquad
\dot p^{(j)}=-h^{(j)}_x(x^{(j)},p^{(j)})
\end{equation*}
subject to the initial condition $\left.(x^{(j)},p^{(j)})\right|_{x^{n+1}=0}=(y,q)$.

\begin{theorem}
\label{theorem about principal symbol of oscillatory integral}
The formula for the principal symbol of the oscillatory integral
$U^{(j)}(x^{n+1})$ reads as follows:
\begin{multline*}
\label{formula for principal symbol of oscillatory integral}
[v^{(j)}(x^{(j)}(x^{n+1};y,q),p^{(j)}(x^{n+1};y,q))]
\,[v^{(j)}(y,q)]^*
\\
\times\exp
\left(
-i\int_0^{x^{n+1}}f^{(j)}(x^{(j)}(\tau;y,q),p^{(j)}(\tau;y,q))\,d\tau
\right).
\end{multline*}
This principal symbol is positively homogeneous in momentum $q$
of degree zero.
\end{theorem}

Theorem \ref{theorem about principal symbol of oscillatory integral}
is due to Safarov \cite{SafarovDSc}. It was later confirmed
by the more detailed analysis carried out in \cite{NicollPhD,jst_part_a}.

Theorem \ref{theorem about principal symbol of oscillatory integral}
is insufficient for the determination of the second term in spectral
asymptotics because one needs information about
the lower order terms of the symbol of the
oscillatory integral $U^{(j)}(x^{n+1})$.
Namely, one needs information about terms
positively homogeneous in momentum $q$
of degree $-1$.
The algorithm described
in Section 2 of \cite{jst_part_a} provides a recursive procedure
for the calculation of all lower order terms, of any degree of homogeneity
in momentum $q$. However, there are two issues here. Firstly, calculations
become very complicated.
Secondly, describing these lower
order terms in an invariant way is problematic.
A few years before his untimely death
Safarov had discussions with one of the authors of this paper and pointed out that,
as far as he was aware, the concept of subprincipal symbol has
never been defined for time-dependent oscillatory integrals
(Fourier integral operators).

We overcome the problem of invariant description of lower order
terms of the symbol of the oscillatory integral $U^{(j)}(x^{n+1})$
by restricting our analysis to $U^{(j)}(0)$. It turns out that
knowing the properties of the lower order terms of the symbol of
$U^{(j)}(0)$ is sufficient for the derivation of two-term spectral asymptotics.
And $U^{(j)}(0)$ is a pseudodifferential operator, so one can use
here the standard notion of subprincipal symbol of a
pseudodifferential operator.

The following result was recently established in \cite{jst_part_a}.

\begin{theorem}
\label{theorem about subprincipal symbol of oscillatory integral}
We have
\begin{equation}
\label{subprincipal symbol of OI at time zero}
\operatorname{tr}[U^{(j)}(0)]_\mathrm{sub}
=-i\{[v^{(j)}]^*,v^{(j)}\},
\end{equation}
where $\,\operatorname{tr}\,$ stands for the matrix trace.
\end{theorem}

It is interesting that the RHS of formula
\eqref{subprincipal symbol of OI at time zero}
admits a geometric interpretation: it can be interpreted
as the scalar curvature of a $\mathrm{U}(1)$ connection
on $T^*M\setminus\{0\}$, see Section 5 of \cite{jst_part_a}
for details. This connection is to do with gauge transformations
of the normalised eigenvector
$v^{(j)}(x,p)$ of the principal
symbol $L_\mathrm{prin}(x,p)$ corresponding to the eigenvalue
$h^{(j)}(x,p)$. Namely, observe that if $v^{(j)}(x,p)$ is an eigenvector
and $\phi^{(j)}(x,p)$ is an arbitrary smooth real-valued function,
then $e^{i\phi^{(j)}(x,p)}v^{(j)}(x,p)$ is also an eigenvector.
Careful analysis of the gauge transformation
\begin{equation}
\label{gauge transformation of the eigenvector}
v^{(j)}\mapsto e^{i\phi^{(j)}}v^{(j)}
\end{equation}
leads to the appearance of
a curvature term.

\section{Mollified spectral asymptotics}
\label{Mollified spectral asymptotics}

Denote by
\begin{equation}
\label{definition of integral kernel of wave group}
u(x,x^{n+1},y):=
\sum_k e^{-ix^{n+1}\lambda_k}v_k(x)[v_k(y)]^*
\end{equation}
the integral kernel of the propagator \eqref{definition of wave group}.
The quantity \eqref{definition of integral kernel of wave group}
can be understood as a matrix-valued distribution in the variable
$x^{n+1}\in\mathbb{R}$ depending on the parameters $x,y\in M$.
Further on in this section
we will be studying the quantity
\begin{equation}
\label{trace of integral kernel of wave group}
\hat f(x,x^{n+1}):=
\operatorname{tr}u(x,x^{n+1},x)=
\sum_k\|v_k(x)\|^2e^{-ix^{n+1}\lambda_k}\,.
\end{equation}
Here $\|v_k(x)\|$ is the Euclidean norm
of the $m$-column $v_k$ evaluated at the point $x\in M$.
Of course, $\|v_k(x)\|^2$ is a real-valued density.

In order to understand the reason for our interest in
\eqref{trace of integral kernel of wave group}, put
\begin{equation}
\label{trace of integral kernel of wave group inverse Fourier}
f(x,\lambda):=
\sum_k\|v_k(x)\|^2\delta(\lambda-\lambda_k).
\end{equation}
Then
\eqref{trace of integral kernel of wave group}
and
\eqref{trace of integral kernel of wave group inverse Fourier}
are related as
$\hat f=\mathcal{F}_{\lambda\to x^{n+1}}[f\,]$
and
$f=\mathcal{F}^{-1}_{x^{n+1}\to\lambda}[\hat f\,]$,
where the one-dimensional Fourier transform $\mathcal{F}$
and its inverse $\mathcal{F}^{-1}$ are defined as in Section 6 of \cite{jst_part_a}.
The quantity \eqref{trace of integral kernel of wave group inverse Fourier}
contains all the information on the spectrum of our operator $L$
and it also contains some information on the eigenfunctions.

Let $\hat\rho:\mathbb{R}\to\mathbb{C}$ be a smooth function such that
$\hat\rho(0)=1$, $\hat\rho'(0)=0$ and the support of $\hat\rho$ is sufficiently
small. Here `sufficiently small' means that
\linebreak
$\operatorname{supp}\hat\rho\subset(-\mathbf{T},\mathbf{T})$,
where $\mathbf{T}$ is the infimum of the lengths of all possible loops.
See Section 6 in \cite{jst_part_a} for details.
Denote also
$\rho(\lambda):=\mathcal{F}^{-1}_{x^{n+1}\to\lambda}[\hat\rho(x^{n+1})]$.

We mollify the distributions
\eqref{trace of integral kernel of wave group}
and
\eqref{trace of integral kernel of wave group inverse Fourier}
by switching to
$\hat\rho(x^{n+1})\hat f(x,x^{n+1})$
and
$(\rho*f)(x,\lambda)$,
where the star indicates convolution in the variable $\lambda$.
It was shown in \cite{jst_part_a}
that Theorems
\ref{theorem about principal symbol of oscillatory integral}
and
\ref{theorem about subprincipal symbol of oscillatory integral}
imply the following result.

\begin{theorem}
\label{mollified theorem 1}
We have, uniformly over $x\in M$,
\begin{equation*}
\label{mollified theorem 1 formula}
(\rho*f)(x,\lambda)=
n\,a(x)\,\lambda^{n-1}+(n-1)\,b(x)\,\lambda^{n-2}+O(\lambda^{n-3})
\quad
\text{as}
\quad
\lambda\to+\infty.
\end{equation*}
Here the densities $a(x)$ and $b(x)$ are given by formulae
\begin{equation}
\label{formula for a}
a(x)=(2\pi)^{-n}\sum_{j=1}^{m/2}
\ \int\limits_{h^{(j)}(x,p)<1}dp\,,
\end{equation}
\begin{multline}
\label{formula for b}
b(x)=-n(2\pi)^{-n}\sum_{j=1}^{m/2}
\ \int\limits_{h^{(j)}(x,p)<1}
\Bigl(
[v^{(j)}]^*L_\mathrm{sub}v^{(j)}
\\
-\frac i2
\{
[v^{(j)}]^*,L_\mathrm{prin}-h^{(j)},v^{(j)}
\}
+\frac i{n-1}h^{(j)}\{[v^{(j)}]^*,v^{(j)}\}
\Bigr)(x,p)\,
dp\,,
\end{multline}
where $dp=dp_1\ldots dp_n$.
\end{theorem}

Theorem \ref{mollified theorem 1} warrants the following remarks.

\begin{remark}
\label{remarks for mollified theorem 1 a}
It is easy to see that
the RHS of formula \eqref{formula for b} is invariant under gauge transformations
\eqref{gauge transformation of the eigenvector} of the eigenvectors of the principal
symbol. 
\end{remark}
\begin{remark}
\label{remarks for mollified theorem 1 b}
Let $R:M\to\mathrm{U}(m)$ be an arbitrary smooth unitary matrix-function.
As one would expect, the RHS of formula \eqref{formula for b}
is invariant under gauge transformations $L\mapsto R^*LR$ of our operator,
but establishing this is not that easy. 
The corresponding calculations are presented in Section 9 of \cite{jst_part_a}.
\end{remark}

Let us now leave in
\eqref{trace of integral kernel of wave group inverse Fourier}
only terms with positive $\lambda_k$ and define the quantity
\begin{equation*}
\label{trace of integral kernel of wave group inverse Fourier plus}
f_+(x,\lambda):=
\sum_{\lambda_k>0}\|v_k(x)\|^2\delta(\lambda-\lambda_k).
\end{equation*}
Theorem \ref{mollified theorem 1} implies the following Corollary.

\begin{corollary}
\label{mollified corollary 1}
We have, uniformly over $x\in M$, the following two results:
\begin{equation*}
\label{mollified corollary 1 formula}
(\rho*f_+)(x,\lambda)=
n\,a(x)\,\lambda^{n-1}+(n-1)\,b(x)\,\lambda^{n-2}+O(\lambda^{n-3})
\quad
\text{as}
\quad
\lambda\to+\infty
\end{equation*}
and $(\rho*f_+)(x,\lambda)$ vanishes faster than any negative power of $|\lambda|$
as $\lambda\to-\infty$.
\end{corollary}

Let us define the two local counting functions
\begin{equation}
\label{definition of local counting functions}
N_\pm(x,\lambda):=
\begin{cases}
0\quad\text{if}\quad\lambda\le0,\\
\sum_{0<\pm\lambda_k<\lambda}\|v_k(x)\|^2\quad\text{if}\quad\lambda>0.
\end{cases}
\end{equation}
The function $N_+(x,\lambda)$ counts the eigenvalues $\lambda_k$
between zero and $\lambda$,
whereas
the function $N_-(x,\lambda)$ counts the eigenvalues $\lambda_k$
between $-\lambda$ and zero.
In both cases counting eigenvalues involves assigning them weights $\|v_k(x)\|^2$.

We have
$(\rho*N_+)(x,\lambda)=\int_{-\infty}^\lambda(\rho*f_+)(x,\mu)\,d\mu$,
so
Corollary~\ref{mollified corollary 1} implies
\begin{multline}
\label{asymptotics for mollified local counting function plus}
(\rho*N_+)(x,\lambda)=
\\
a(x)\,\lambda^n+b(x)\,\lambda^{n-1}+
\begin{cases}
O(\lambda^{n-2})&\text{if}\quad n\ge3,\\
O(\ln\lambda)&\text{if}\quad n=2,
\end{cases}
\quad
\text{as}
\quad
\lambda\to+\infty.
\end{multline}

The asymptotics for $(\rho*N_-)(x,\lambda)$ is obtained by applying the above result to
the operator $-L$ and using the symmetries \eqref{symmetry property 1}. This gives
\begin{multline}
\label{asymptotics for mollified local counting function minus}
(\rho*N_-)(x,\lambda)=
\\
a(x)\,\lambda^n-b(x)\,\lambda^{n-1}+
\begin{cases}
O(\lambda^{n-2})&\text{if}\quad n\ge3,\\
O(\ln\lambda)&\text{if}\quad n=2,
\end{cases}
\quad
\text{as}
\quad
\lambda\to+\infty.
\end{multline}

Note that the second terms in the asymptotic formulae
\eqref{asymptotics for mollified local counting function plus}
and
\eqref{asymptotics for mollified local counting function minus}
have opposite signs
and that the remainders are uniform in $x\in M$.

Finally, let us define the two global counting functions
\begin{equation}
\label{definition of global counting functions}
N_\pm(\lambda):=
\begin{cases}
0\quad\text{if}\quad\lambda\le0,\\
\sum_{0<\pm\lambda_k<\lambda}1\quad\text{if}\quad\lambda>0.
\end{cases}
\end{equation}
We have $N_\pm(\lambda)=\int_M N_\pm(x,\lambda)\,dx$, where $dx=dx^1\ldots dx^n$.
Therefore, formulae
\eqref{asymptotics for mollified local counting function plus}
and
\eqref{asymptotics for mollified local counting function minus}
imply
\begin{equation}
\label{asymptotics for mollified global counting function plus minus}
(\rho*N_\pm)(\lambda)=
a\,\lambda^n\pm b\,\lambda^{n-1}+
\begin{cases}
O(\lambda^{n-2})&\text{if}\quad n\ge3,\\
O(\ln\lambda)&\text{if}\quad n=2,
\end{cases}
\quad
\text{as}
\quad
\lambda\to+\infty,
\end{equation}
where
\begin{equation}
\label{integration of local Weyl coefficients}
a=\int_Ma(x)\,dx\,,
\qquad
b=\int_Mb(x)\,dx\,.
\end{equation}

\section{Unmollified spectral asymptotics}
\label{Unmollified spectral asymptotics}

In this section we write down asymptotic formulae for the local and global counting functions without mollification.
These can be obtained from the mollified asymptotic formulae
\eqref{asymptotics for mollified local counting function plus},
\eqref{asymptotics for mollified local counting function minus}
and
\eqref{asymptotics for mollified global counting function plus minus}
by applying appropriate Fourier Tauberian theorems,
see Section 8 of \cite{jst_part_a}.

\begin{theorem}
\label{theorem spectral function unmollified one term}
We have, uniformly over $x\in M$,
\begin{equation}
\label{theorem spectral function unmollified one term formula}
N_\pm(x,\lambda)=a(x)\,\lambda^n+O(\lambda^{n-1})
\quad
\text{as}
\quad
\lambda\to+\infty.
\end{equation}
\end{theorem}

\begin{corollary}
\label{theorem counting function unmollified one term}
We have
\begin{equation}
\label{theorem counting function unmollified one term formula}
N_\pm(\lambda)=a\lambda^n+O(\lambda^{n-1})
\quad
\text{as}
\quad
\lambda\to+\infty.
\end{equation}
\end{corollary}

\begin{theorem}
\label{theorem spectral function unmollified two term}
If the point $x\in M$ is nonfocal then
\begin{equation}
\label{two-term asymptotic formula for spectral function}
N_\pm(x,\lambda)=a(x)\,\lambda^n\pm b(x)\,\lambda^{n-1}+o(\lambda^{n-1})
\quad
\text{as}
\quad
\lambda\to+\infty.
\end{equation}
\end{theorem}

\begin{theorem}
\label{theorem counting function unmollified two term}
If the nonperiodicity condition is fulfilled then
\begin{equation}
\label{two-term asymptotic formula for counting function}
N_\pm(\lambda)=a\lambda^n\pm b\lambda^{n-1}+o(\lambda^{n-1})
\quad
\text{as}
\quad
\lambda\to+\infty.
\end{equation}
\end{theorem}

A point $x\in M$ is said to be nonfocal if there are not too many
Hamiltonian loops originating from this point. `Nonperiodicity' means that
there are not too many periodic Hamiltonian trajectories. See subsection 8.2
in \cite{jst_part_a} for details.

The results presented in this section were first obtained by Victor Ivrii
\cite{ivrii_springer_lecture_notes,ivrii_book}
but without an explicit formula for the second asymptotic coefficient.

Asymptotic formulae of the type
\eqref{theorem spectral function unmollified one term formula}--\eqref{two-term asymptotic formula for counting function}
are called \emph{Weyl-type formulae} and the
coefficients in such formulae are often
referred to as \emph{Weyl coefficients}.

\section{The eta function}
\label{The eta function}

The (global) eta function of our operator $L$ is defined as
\begin{equation}
\label{definition of eta function}
\eta(s):=\sum_{\lambda_k\ne0}
\frac{\operatorname{sgn}\lambda_k}{|\lambda_k|^s}
=\int_0^{+\infty}\lambda^{-s}\,(N_+'(\lambda)-N_-'(\lambda))\,d\lambda\,,
\end{equation}
where summation is carried out over all nonzero eigenvalues $\lambda_k$ of $L$
and $s\in\mathbb{C}$ is the independent variable.
The series \eqref{definition of eta function}
converges absolutely for $\operatorname{Re}s>n$
and defines a holomorphic function in
this half-plane. Moreover, it is known~\cite{atiyah_part_3}
that it
extends meromorphically
to the whole $s$-plane with simple poles which can only occur at real integer
values of $s$.
The eta function is the accepted way of describing the asymmetry of the spectrum.

Formula \eqref{theorem counting function unmollified one term formula}
implies that the eta function does not have a pole at $s=n$
and formula
\eqref{asymptotics for mollified global counting function plus minus}
implies that at $s=n-1$ the residue is
\begin{equation}
\label{residue at n minus 1}
\operatorname{Res}(\eta\,,n-1)=2(n-1)b\,,
\end{equation}
where $b$ is the coefficient from
\eqref{asymptotics for mollified global counting function plus minus}.
Thus, for a generic operator $L$ the first pole of the eta function
is at $s=n-1$ and formulae
\eqref{formula for b}
\eqref{integration of local Weyl coefficients}
and
\eqref{residue at n minus 1}
give us an explicit expression for the residue.

It is known \cite{atiyah_part_3,gilkey_1981}
that the eta function does not have a pole at $s=0$.
The real number $\eta(0)$ is called the \emph{eta invariant}.
It can be interpreted as
the number of positive eigenvalues minus the number of negative eigenvalues.
This interpretation is based on the observation that if we
were dealing with an Hermitian matrix $L$,
then $\eta(0)$ would indeed be
the number of positive eigenvalues minus the number of negative eigenvalues.
Our differential operator $L$ has infinitely many positive eigenvalues and
infinitely many negative eigenvalues, and the concept of the eta function allows us to
regularise the expression
`number of positive eigenvalues minus the number of negative eigenvalues'.

The eta function may have poles at
\begin{equation}
\label{admissible1}
s=n-1,\ldots,1,-1,-2,\ldots.
\end{equation}
However, a more careful analysis \cite{romanian}
shows that poles may occur only at values of $s$ of the form
\begin{equation}
\label{admissible2}
s=n-1-2l,\qquad l=0,1,\ldots.
\end{equation}
The authors of \cite{romanian} call values of $s$ from the intersection of the sets
\eqref{admissible1} and \eqref{admissible2} \emph{admissible}.
It was shown in
\cite{romanian}
that residues of the eta function at positive admissible integers
are generically nonzero. This agrees with our explicit calculation
of the residue at $s=n-1$.

\section{Systems of two equations}
\label{Systems of two equations}

From now on we assume that
\begin{equation}
\label{m equals two}
m=2
\end{equation}
and that
\begin{equation}
\label{principal symbol is trace-free}
\operatorname{tr}L_\mathrm{prin}(x,p)=0.
\end{equation}
In other words, we now restrict our analysis to $2\times2$ operators
with trace-free principal symbols. The logic behind the assumptions
\eqref{m equals two} and \eqref{principal symbol is trace-free}
is that they single out the simplest class of first order systems
and we expect to extract more geometry out of our asymptotic analysis
and simplify the results.

It is easy to see that formulae
\eqref{definition of ellipticity},
\eqref{m equals two} and \eqref{principal symbol is trace-free}
imply that the dimension $n$ of our manifold $M$ is less than or equal to three.
Further on we assume that
\begin{equation}
\label{n equals three}
n=3.
\end{equation}

\begin{remark}
\label{remark on  parallelizability}
It was shown in \cite{jst_part_b}
that a 3-manifold is parallelizable if and only if
there exists a self-adjoint elliptic first order linear differential operator
with trace-free principal symbol acting on 2-columns of
complex-valued half-densities over this manifold.
This means that once we restricted our analysis to the special case
\eqref{m equals two}--\eqref{n equals three} we are working on
a parallelizable manifold.
\end{remark}

\begin{remark}
\label{remark on orientability}
It is well known that a 3-manifold is orientable if and only if it is parallelizable,
see Theorem 1 in Chapter VII of \cite{Kirby}.
\end{remark}

Observe that under the assumption \eqref{m equals two}
the determinant of the principal symbol is a quadratic form
in the dual variable (momentum) $p$\,:
\begin{equation}
\label{definition of metric}
\det L_\mathrm{prin}(x,p)=-g^{\alpha\beta}(x)\,p_\alpha p_\beta\,.
\end{equation}
Furthermore, formulae
\eqref{definition of ellipticity}
and
\eqref{principal symbol is trace-free}
imply that the quadratic form
$g^{\alpha\beta}(x)\,p_\alpha p_\beta\,$
is positive definite.
We interpret the real coefficients $g^{\alpha\beta}(x)=g^{\beta\alpha}(x)$,
$\alpha,\beta=1,2,3$,
as components of a (contravariant) metric tensor.
Thus, $2\times2$ operators with trace-free principal symbols
are special in that the concept of a Riemannian metric is encoded
within such operators. This opens the way to the geometric
interpretation of our analytic results.

Note also that under the assumptions
\eqref{m equals two} and \eqref{principal symbol is trace-free}
the principal symbol of the operator $L^2$
is automatically proportional to the $2\times2$ identity
matrix $I$:
\begin{equation}
\label{Dirac-type operators}
(L^2)_\mathrm{prin}(x,p)
=(L_\mathrm{prin})^2(x,p)
=\bigl(g^{\alpha\beta}(x)\,p_\alpha p_\beta\bigr)\,I\,.
\end{equation}
Operators possessing the property
\eqref{Dirac-type operators}
are called \emph{Dirac-type operators}.

Now take an arbitrary smooth matrix-function
\begin{equation}
\label{SU2 matrix-function R}
R:M\to\mathrm{SU}(2)
\end{equation}
and consider the transformation of our $2\times2$ differential operator
\begin{equation}
\label{SU2 transformation of the operator}
L\mapsto R^*LR\,.
\end{equation}
We interpret
\eqref{SU2 transformation of the operator} as a gauge transformation
because it does not affect our counting functions
\eqref{definition of local counting functions},
\eqref{definition of global counting functions}
and the eta function \eqref{definition of eta function}.
Note also that the transformation \eqref{SU2 transformation of the operator}
preserves the condition \eqref{principal symbol is trace-free}.

The transformation \eqref{SU2 transformation of the operator}
of the differential operator $L$ induces the following transformations
of its principal and subprincipal symbols:
\begin{equation}
\label{SU2 transformation of the principal symbol}
L_\mathrm{prin}\mapsto R^*L_\mathrm{prin}R\,,
\end{equation}
\begin{equation}
\label{SU2 transformation of the subprincipal symbol}
L_\mathrm{sub}\mapsto
R^*L_\mathrm{sub}R
+\frac i2
\left(
R^*_{x^\alpha}(L_\mathrm{prin})_{p_\alpha}R
-
R^*(L_\mathrm{prin})_{p_\alpha}R_{x^\alpha}
\right).
\end{equation}
Comparing formulae
\eqref{SU2 transformation of the principal symbol}
and
\eqref{SU2 transformation of the subprincipal symbol}
we see that, unlike the principal symbol, the subprincipal
symbol does not transform in a covariant fashion due to
the appearance of terms with the gradient of the
matrix-function $R(x)$. 

It turns out that one can overcome the non-covariance
in \eqref{SU2 transformation of the subprincipal symbol} by introducing
the \emph{covariant subprincipal symbol} $\,L_\mathrm{csub}(x)\,$
in accordance with formula
\begin{equation}
\label{definition of covariant subprincipal symbol}
L_\mathrm{csub}:=
L_\mathrm{sub}
-\frac i{16}\,
g_{\alpha\beta}
\{
L_\mathrm{prin}
,
L_\mathrm{prin}
,
L_\mathrm{prin}
\}_{p_\alpha p_\beta},
\end{equation}
where subscripts $p_\alpha$ and $p_\beta$ indicate partial derivatives and
curly brackets denote the generalised Poisson bracket on matrix-functions
\eqref{definition of generalised Poisson bracket}.

\begin{lemma}
\label{Lemma about covariant subprincipal symbol}
The transformation
\eqref{SU2 transformation of the operator}
of the differential operator induces the transformation
$\,L_\mathrm{csub}\mapsto R^*L_\mathrm{csub}R\,$
of its covariant subprincipal symbol.
\end{lemma}

The proof of
Lemma~\ref{Lemma about covariant subprincipal symbol}
was given in \cite{nongeometric}.
Note that the analysis in \cite{nongeometric} was performed in a more general,
4-dimensional Lorentzian setting.

In our 3-dimensional Riemannian setting the correction term in the RHS of
\eqref{definition of covariant subprincipal symbol} turns out to be proportional to the
$2\times2$ identity matrix $I$:
\begin{equation}
\label{correction term proportional to the identity matrix}
-\frac i{16}\,
g_{\alpha\beta}
\{
L_\mathrm{prin}
,
L_\mathrm{prin}
,
L_\mathrm{prin}
\}_{p_\alpha p_\beta}=If\,,
\end{equation}
where $f:M\to\mathbb{R}$ is some scalar function.
The function $f(x)$ in \eqref{correction term proportional to the identity matrix}
admits a geometric interpretation \cite{jst_part_b,action}
but we do not discuss this in the current paper.

\begin{theorem}
\label{mollified theorem 1 special}
In the special case
\eqref{m equals two}--\eqref{n equals three}
formulae
\eqref{formula for a}
and
\eqref{formula for b}
read
\begin{equation}
\label{formula for a special}
a(x)=\frac1{6\pi^2}\,
\sqrt{\det g_{\alpha\beta}(x)}\,,
\end{equation}
\begin{equation}
\label{formula for b special}
b(x)=-\frac1{4\pi^2}\,
\bigl(
(\operatorname{tr}L_\mathrm{csub})\,\sqrt{\det g_{\alpha\beta}}\,
\bigr)(x)\,.
\end{equation}
\end{theorem}

Theorem \eqref{mollified theorem 1 special} was
established in \cite{jst_part_b}, though the density $b(x)$
was written in \cite{jst_part_b} in a slightly different way.
The use of the concept of covariant subprincipal symbol introduced
in \cite{nongeometric} allows us to replace
formula (1.19) from \cite{jst_part_b} by the more compact expression
\eqref{formula for b special}.

Formula \eqref{formula for a special} tells us that the first local Weyl coefficient
is proportional to the standard Riemannian density. The first global Weyl coefficient
is obtained from the local one by integration, see formula
\eqref{integration of local Weyl coefficients}, and is proportional to
the Riemannian volume $V$ of our manifold $M$:
$\,a=\frac1{6\pi^2}\,V\,$.

In order to understand the geometric meaning of formula
\eqref{formula for b special} we observe that the covariant subprincipal
symbol can be uniquely represented in the form
\begin{equation}
\label{decomposition of  covariant subprincipal symbol}
L_\mathrm{csub}(x)=L_\mathrm{prin}(x,A(x))+IA_4(x),
\end{equation}
where $A=(A_1,A_2,A_3)$ is some real-valued covector field,
$A_4$ is some real-valued scalar field, $x=(x^1,x^2,x^3)$ are local coordinates on $M$
(we are working in the nonrelativistic setting) and $I$ is the $2\times2$ identity matrix.
Applying the results of \cite{nongeometric} to the relativistic operator
appearing in the LHS of \eqref{dynamic equation most basic}
we conclude that $A=(A_1,A_2,A_3)$ is the magnetic covector potential and
$A_4$ is the electric potential.
Note that
Lemma~\ref{Lemma about covariant subprincipal symbol}
and formulae
\eqref{SU2 transformation of the principal symbol}
and
\eqref{decomposition of  covariant subprincipal symbol}
tell us that the magnetic covector potential and electric potential
are invariant under gauge transformations
\eqref{SU2 transformation of the operator}.

Substituting \eqref{decomposition of  covariant subprincipal symbol}
into \eqref{formula for b special}
and making use of \eqref{principal symbol is trace-free}
we get
\begin{equation}
\label{formula for b special in terms of electric potential}
b(x)=-\frac1{2\pi^2}\,
\bigl(
A_4\,\sqrt{\det g_{\alpha\beta}}\,
\bigr)(x)\,.
\end{equation}
Thus, the second Weyl coefficient is proportional to the electric potential
and does not depend on the magnetic covector potential.

A number of researchers have studied the effect of the electromagnetic field
on the spectrum of the first order differential operator $L$ under the assumptions
\eqref{m equals two}--\eqref{n equals three} and our formula
\eqref{formula for b special in terms of electric potential}
is a further contribution to this line of research.
However, we believe that such results do not have a physical meaning because
our $2\times2$ first order differential operator $L$ describes a massless particle
and no known massless particle has an electric charge.
In the absence of an electric charge
the particle cannot interact with the electromagnetic field.

The electron is an example of a charged massive
particle but it is described by a $4\times4$ first order differential operator.
Also, in the case of the electron it is more natural to do asymptotic analysis
in a different setting, with Planck's constant tending to zero.
Spectral problems for the electron in 3-dimensional Euclidean space
in the presence of magnetic and electric potentials
were extensively
studied by Ivrii \cite{ivrii_springer_lecture_notes,ivrii_book}.
An analytic
(i.e.~based on the concepts of principal symbol and covariant subprincipal symbol)
representation of the massive Dirac equation in curved 4-dimensional
Lorentzian spacetime was given in \cite{nongeometric}.

\section{Spin structure}
\label{Spin structure}

Let $M$ be a connected closed oriented Riemannian 3-manifold.
Let us consider all possible self-adjoint elliptic first order $2\times2$
linear differential operators $L$ with trace-free principal symbols
corresponding, in the sense of formula \eqref{definition of metric},
to the prescribed metric.
See also Remarks
\ref{remark on  parallelizability}
and
\ref{remark on orientability}.
In this section our aim is to classify all such operators $L$.

We define the topological charge as
\begin{equation}
\label{definition of relative orientation}
\mathbf{c}:=-\frac i2\sqrt{\det g_{\alpha\beta}}\,\operatorname{tr}
\bigl((L_\mathrm{prin})_{p_1}(L_\mathrm{prin})_{p_2}(L_\mathrm{prin})_{p_3}\bigr),
\end{equation}
with the subscripts $p_1$, $p_2$, $p_3$
indicating partial derivatives with respect to the components of momentum
$p=(p_1,p_2,p_3)$.
It was shown in Section 3 of \cite{jst_part_b} that
the number $\mathbf{c}$ defined by formula
(\ref{definition of relative orientation})
can take only two values, $+1$ or $-1$,
and describes the orientation of the principal symbol
relative to the chosen orientation of local coordinates $x=(x^1,x^2,x^3)$.
Of course, the transformation $L\mapsto-L$ inverts the topological charge.

Further on we work with operators whose topological charge is $+1$.

We say that the operators $L$ and $\tilde L$ are equivalent if there
exists a smooth matrix-function \eqref{SU2 matrix-function R}
such that $\tilde L_\mathrm{prin}=R^*L_\mathrm{prin}R$.
The equivalence classes of operators obtained this way
are called \emph{spin structures}.

An example of non-equivalent operators $L$ and $\tilde L$ on the 3-torus
was given in Appendix~A of \cite{jst_part_b}. Furthermore, using the above
definition of spin structure one can show that there are eight
distinct spin structures on the 3-torus whereas the spin structure
on the 3-sphere is unique.

We see that an operator $L$ is uniquely determined,
modulo a gauge transformation
\eqref{SU2 transformation of the operator},
by the metric, topological charge, spin structure,
magnetic covector potential and electric potential.

We claim that in dimension three our analytic definition of spin structure
is equivalent to the traditional topological definition.
We will provide a rigorous proof of this claim in a separate paper.

\section{The massless Dirac operator}
\label{The massless Dirac operator}

In this section we continue dealing with the special case
\eqref{m equals two}--\eqref{n equals three}
but make the additional assumption that the
magnetic covector potential and electric potential vanish.
The resulting operator $L$ is called the \emph{massless Dirac operator on half-densities}.
It is uniquely determined,
modulo a gauge transformation
\eqref{SU2 transformation of the operator},
by the metric, topological charge and spin structure.

In practice most researchers work with the
massless Dirac operator which acts on 2-columns of
complex-valued scalar fields (components of a Weyl spinor)
rather than on 2-columns of complex-valued half-densities.
As we have a Riemannian metric encoded in the principal symbol
of our operator, scalar fields can be identified with half-densities:
it is just a matter of multiplying or dividing by
$(\det g_{\alpha\beta})^{1/4}$.
Hence, the `traditional' massless Dirac operator and
the massless Dirac operator on half-densities are related by a
simple formula, see formula (A.19) in \cite{jst_part_b},
and their spectra are the same.
For spectral theoretic purposes it is more convenient to work
with half-densities because in this case the inner product does not depend
on the metric.

The massless Dirac operator describes the massless neutrino.
We are looking at a single neutrino living in a closed 3-dimensional
Riemannian universe. The eigenvalues are the energy levels of the particle.
The tradition is to associate positive eigenvalues with the energy levels of the neutrino
and negative eigenvalues with the energy levels  of the antineutrino.

Formula \eqref{formula for b special in terms of electric potential}
tells us that the second Weyl coefficient for the massless Dirac operator
is zero, both locally and globally.
Formula \eqref{residue at n minus 1},
in turn, tells us that the eta function of the massless Dirac operator
does not have a pole at $s=2$.

The natural question is where is the first pole of the eta function?
It was shown in \cite{bismut_and_freed}
that the eta function of the massless Dirac operator is holomorphic
in the half-plane $\operatorname{Re}s>-2$.
This agrees with formulae
\eqref{admissible1}
and
\eqref{admissible2}.

Furthermore, Branson and Gilkey \cite{branson_and_gilkey} have shown
that generically the eta function of the massless Dirac operator
has a pole at $s=-2$ and calculated the residue. Consider the covariant rank three
tensor $(\nabla_{\alpha}Ric_{\beta\nu})Ric_\gamma{}^\nu$,
where $\nabla$ stands for the covariant derivative and $Ric$ for Ricci curvature
(both are understood in terms of the Levi-Civita connection),
and antisymmetrize it. This gives a totally antisymmetric covariant
rank three tensor which is equivalent to a 3-form. According to
\cite{branson_and_gilkey}, the integral of this 3-form over the 3-manifold $M$
gives, up to a particular nonzero constant factor,
the residue of the eta function of the massless Dirac operator
at $s=-2$.

The fact that the first pole of the eta function of the massless Dirac operator
is at $s=-2$ indicates that with a very high accuracy
the large (in terms of modulus) positive and negative eigenvalues
are distributed in the same way.
This, in turn, means that the massless Dirac operator is special
and has hidden symmetries encoded in it.

We end this section by highlighting one particular symmetry of the
massless Dirac operator. Consider the following antilinear operator
acting on 2-columns of complex-valued half-densities:
\begin{equation}
\label{antilinear}
v=
\begin{pmatrix}
v_1\\
v_2
\end{pmatrix}
\mapsto
\begin{pmatrix}
-\overline{v_2}\\
\overline{v_1}
\end{pmatrix}
=:\mathrm{C}(v).
\end{equation}
The operator $\mathrm{C}$ defined by formula \eqref{antilinear}
is called the \emph{charge conjugation operator}.
It is known, see Appendix A in \cite{jst_part_b},
that the linear massless Dirac operator on half-densities $L$ and
the antilinear charge conjugation operator $\mathrm{C}$ commute:
\begin{equation}
\label{commutes}
\mathrm{C}(Lv)=L\mathrm{C}(v).
\end{equation}
Formula \eqref{commutes}
implies, in particular,
that all eigenvalues of the massless Dirac operator have even multiplicity.

The addition of an electric potential preserves the symmetry \eqref{commutes},
but the addition of a magnetic covector potential destroys it.

\section{Small eigenvalues}
\label{Small eigenvalues}

Up till now we dealt with large, in terms of modulus, eigenvalues.
In this section we will deal with small eigenvalues of the
massless Dirac operator.

Suppose that we are working on a connected closed oriented Riemannian 3-manifold
and let $\lambda^{(0)}$ be a double eigenvalue of the massless Dirac operator.
As explained in the end of the previous section, multiplicity two is
the lowest possible.
We now perturb the metric, i.e.~consider an arbitrary metric
$g_{\alpha\beta}(x;\epsilon)$ the components of which are smooth
functions of local coordinates $x^\alpha$,
$\alpha=1,2,3$, and small real parameter $\epsilon$;
here we assume that for $\epsilon=0$ we get the original metric.
In this case one can expand the eigenvalue into an asymptotic series
in powers of the small parameter $\epsilon$:
$\,\lambda(\epsilon)=\lambda^{(0)}
+\lambda^{(1)}\epsilon
+\lambda^{(2)}\epsilon^2
+\ldots\,$
with some constants $\lambda^{(1)},\lambda^{(2)},\ldots\,$.
This asymptotic construction was described
in Sections 3--5 of \cite{torus}. The construction
is somewhat nontrivial because we are dealing with
a double eigenvalue that cannot split.

We now consider two special cases.
In both cases the unperturbed spectrum is symmetric
but symmetry is broken under generic perturbations of the metric.

\subsection{The 3-torus with standard spin structure}

Here the unperturbed metric is assumed to be Euclidean
and standard spin structure means that our unperturbed
massless Dirac operator can be written as an operator with
constant coefficients in the natural $2\pi$-periodic
cyclic coordinates parameterizing the 3-torus,
see formula (1.1) in \cite{torus}.

The spectrum of the unperturbed operator is known,
see, for example, Appendix B in \cite{jst_part_b}
or Section 1 in \cite{torus}.
The smallest eigenvalue is the double eigenvalue $\lambda^{(0)}=0$.
It was shown in \cite{torus} that
\begin{equation}
\label{result for torus}
\lambda(\epsilon)=
\lambda^{(2)}\epsilon^2
+O(\epsilon^3)
\quad\text{as}\quad\epsilon\to0
\end{equation}
with an explicit expression for the constant $\lambda^{(2)}$.
Examination of this explicit expression shows
that under a generic perturbation of the metric
we get
$\,\lambda^{(2)}\ne0\,$
which is an indication of spectral asymmetry.

Furthermore, two special families of metrics were identified in \cite{torus}
for which the eigenvalue closest to zero,
$\lambda(\epsilon)$, can be evaluated explicitly.
Formula \eqref{result for torus}
was tested against explicit results for these two families of metrics.

\subsection{The 3-sphere}

Here the unperturbed metric is obtained by restricting the Euclidean
metric from $\mathbb{R}^4$ to $\mathbb{S}^3$.
There is no issue with spin structure because for the 3-sphere the
spin structure is unique.

The spectrum of the unperturbed operator is known,
see, for example, Appendix B in \cite{jst_part_b}.
The smallest, in terms of modulus, eigenvalues are
the double eigenvalues $\lambda_+^{(0)}=+\frac32$ and $\lambda_-^{(0)}=-\frac32\,$.
We get
\begin{equation}
\label{result for sphere}
\lambda_\pm(\epsilon)=
\pm\frac32
+\lambda_\pm^{(1)}\epsilon
+\lambda_\pm^{(2)}\epsilon^2
+O(\epsilon^3)
\quad\text{as}\quad\epsilon\to0.
\end{equation}

In order to write down the coefficients $\lambda_\pm^{(1)}$
we consider the Riemannian volume $V(\epsilon)$ of our manifold $M$
and expand it in powers of $\epsilon$:
\begin{equation}
\label{asymptotic expansion for volume}
V(\epsilon)=V^{(0)}+V^{(1)}\epsilon
+O(\epsilon^2)
\qquad\text{as}\qquad\epsilon\to0,
\end{equation}
where
$V^{(0)}=\,2\pi^2$
is the volume of the unperturbed 3-sphere.
It turns out that
\begin{equation}
\label{formula for lambda1}
\lambda_\pm^{(1)}=\mp\frac1{4\pi^2}V^{(1)}.
\end{equation}
Formulae
\eqref{result for sphere}--\eqref{formula for lambda1}
tell us that in the first approximation
in $\epsilon$ spectral symmetry is preserved and
the increments of the two eigenvalues
closest to zero are determined
by the increment of volume.
If the volume increases then the moduli of the two eigenvalues
closest to zero decrease and in the first approximation in
$\epsilon$ they decrease in the same way.

Arguing along the lines of \cite{torus} one can write down
explicit expressions for the constants $\lambda_\pm^{(2)}$.
Examination of these explicit expressions shows
that under a generic perturbation of the metric
we get spectral asymmetry in the $\epsilon^2$ terms:
$\,\lambda_+^{(2)}+\lambda_-^{(2)}\ne0\,$.
A detailed exposition will be provided in a separate paper.

Note that there is a family of metrics for which
the two eigenvalues closest to zero,
$\lambda_+(\epsilon)$ and $\lambda_-(\epsilon)$,
can be evaluated explicitly.
These are the so-called generalized Berger metrics: see
Proposition 3.1 in \cite{hitchin}
or Definition~4 in~\cite{godbout}.

\section{Issue with eigenvalues of the principal symbol}
\label{Issue with eigenvalues of the principal symbol}

Throughout this paper we assumed that
the eigenvalues of the matrix-function $L_\mathrm{prin}(x,p)$,
the principal symbol of our operator $L$,
are simple for all $(x,p)\in T^*M\setminus\{0\}$.
In this section we briefly examine the issues that arise if one drops this assumption.

Ivrii showed that Theorem~\ref{theorem spectral function unmollified one term}
holds without any assumptions on the eigenvalues of the principal symbol,
see
Theorem 0.1 in \cite{IvriiFuncAn1982}
or 
Theorem 0.1 in \cite{ivrii_springer_lecture_notes}.
However, establishing analogues of
Theorems
\ref{theorem spectral function unmollified two term}
and
\ref{theorem counting function unmollified two term}
without the assumption that the eigenvalues of
the principal symbol are simple is not, by any means, straightforward and there
are two issues that have to be addressed.
These are highlighted in the following two subsections.

\subsection{Geometric conditions for the existence of two-term spectral asymptotics}

If the multiplicity of eigenvalues of the principal symbol varies as a function
of $(x,p)\in T^*M\setminus\{0\}$ then the expectation is that one needs
to consider `generalised' Hamiltonian trajectories, with branching occurring
at points in the cotangent bundle
where multiplicities of eigenvalues of the principal symbol change.

Ivrii \cite{IvriiFuncAn1982,ivrii_springer_lecture_notes}
dealt with the issue of variable multiplicities of eigenvalues of the
principal symbol by assuming that the set of Hamiltonian trajectories
encountering points
where multiplicities of eigenvalues of the principal symbol change
is, in some sense, small.

G.V.~Rozenblyum \cite{grisha}
and later
I.~Kamotski and M.~Ruzhansky \cite{kamotski}
considered `generalised' Hamiltonian trajectories with branching
assuming that the principal
symbol of the operator is well behaved
at points
where multiplicities of eigenvalues of the principal symbol change.
Here good behaviour is understood as smooth microlocal diagonalisability
of the principal symbol plus some conditions on the Poisson brackets
of eigenvalues.

\subsection{Explicit formulae for the second Weyl coefficient}

In the case when the eigenvalues of the principal symbol
are not simple
explicit formulae for the second Weyl coefficient are not known.

A good starting point for the derivation of such formulae would be the
analysis of the case when eigenvalues of the principal symbol
have constant multiplicities for all $(x,p)\in T^*M\setminus\{0\}$.
Let $L_\mathrm{prin}(x,p)$,
be our $m\times m$ principal symbol
and let
$l_j$, $j=1,\ldots,k$, be the multiplicities of its positive eigenvalues,
so that $l_1+\ldots+l_k=m/2$.
Then one can, by analogy with Section 5 of \cite{jst_part_a},
introduce a $\mathrm{U}(l_j)$ connection associated with the $j$th
positive eigenspace of the principal symbol.
It is natural to conjecture that the curvature of this
$\mathrm{U}(l_j)$ connection will appear in the explicit formula for
the second Weyl coefficient.

\section*{Acknowledgments}

The authors are grateful to Nikolai Saveliev for helpful advice.


\end{document}